\documentclass[reqno,centertags,11 pt]{amsart}
\usepackage{amsmath,amsthm,amscd}
\usepackage{amssymb,latexsym}
\usepackage{diagrams}
\usepackage{ucs}
\usepackage[dvips]{epsfig}
\usepackage{array}
\newcommand{\F}{\mathbb F}
\newcommand{\Fpbar}{\overline{\F}_p}
\newcommand{\R}{\mathbb R}
\newcommand{\C}{\mathbb C} 
\newcommand{\Q}{\mathbb Q}
\newcommand{\Z}{\mathbb Z}  

\newcommand{\A}{\mathbb A}

\newcommand{\epsi}{\varepsilon}

\newcommand{\Res}{\mathrm{Res}}

\newcommand{\Frob}{\mathrm{Frob}}

\newcommand{\Pic}{\mathrm{Pic}}
\newcommand{\Gal}{\mathrm{Gal}}

\newcommand{\Sym}{\mathrm{Sym}}

\newcommand{\Spec}{\mathrm{Spec\mbox{ } }}
\newcommand{\GL}{\mathrm{GL}}
\newcommand{\p}{\mathfrak{p}}

\newcommand{\tensor}{\otimes}
\newtheorem{theorem}{Theorem}[section]
\newtheorem{proposition}[theorem]{Proposition}
\newtheorem{lemma}[theorem]{Lemma}

\theoremstyle{definition}
\newtheorem{example}{Example}
\newtheorem{conjecture}[example]{Conjecture}
\newtheorem{definition}[theorem]{Definition}

\theoremstyle{remark}
\newtheorem{remark}[theorem]{Remark}

\begin{document} 
\title[Weights in Serre's conjecture for Hilbert modular forms]{Weights in Serre's conjecture for Hilbert modular forms: the ramified case}
\author[Michael M. Schein]{Michael M. Schein}
\address{Institute of Mathematics, Hebrew University of Jerusalem, Givat Ram, Jerusalem 91904, ISRAEL}
\email{mschein@math.huji.ac.il}
\thanks{The author thanks the NSF for a Graduate Research Fellowship that supported him during part of this work.}
\date{October 16, 2006; revised February 14, 2007}
\begin{abstract}
Let $F$ be a totally real field and $p \geq 3$ a prime.  If $\rho: \Gal(\overline{F}/F) \to \GL_2(\Fpbar)$ is continuous, semisimple, totally odd, and tamely ramified at all places of $F$ dividing $p$, then we formulate a conjecture specifying the weights for which $\rho$ is modular.  This extends the conjecture of Diamond, Buzzard, and Jarvis, which required $p$ to be unramified in $F$.  We also prove a theorem that verifies one half of the conjecture in many cases and use Demb\'{e}l\'{e}'s computations of Hilbert modular forms over $\Q(\sqrt{5})$ to provide evidence in support of the conjecture.
\end{abstract}

\maketitle
\section{Introduction}
Let $F$ be a totally real field and $p \geq 3$ a rational prime.  For any place $v$ of $F$, we write $\mathcal{O}_v$ for the completion of $\mathcal{O}_F$ at $v$ and $k_v$ for the residue field.  Let $p \mathcal{O}_F = \prod_{v | p} v^{e_v}$ be the factorization of $p$ into prime ideals of $F$, so that $e_v$ is the ramification index of $F_v$ over $\Q_p$.  The purpose of this paper is to formulate, and prove some cases of, a Serre-type ``epsilon conjecture'' for mod $p$ Hilbert modular forms over $F$.  Previously this has been done only in the case of $p$ unramified in $F$, i.e. $e_v = 1$ for all $v | p$.
\begin{definition}
A {\emph {(Serre) weight}} is an irreducible $\Fpbar$-representation of the group $\GL_2(\mathcal{O}_F /p) = \prod_{v |p} \GL_2(\mathcal{O}_F / v^{e_v})$.
\end{definition}

Any irreducible mod $p$ representation of $\GL_2(\mathcal{O}_F / v^{e_v})$ factors through the natural surjection $\GL_2(\mathcal{O}_F /v^{e_v}) \to \GL_2(k_v)$; indeed, the kernel is a $p$-group and hence acts trivially (see \cite{EdixhovenFLT} for a proof of this).  By Proposition 1 of \cite{BarthelLivne}, the irreducible $\Fpbar$-representations of $\GL_2(k_v)$ are:
$$ \sigma_v = \bigotimes_{\tau: k_v \hookrightarrow \Fpbar} (\det\nolimits^{w_\tau} \Sym^{k_\tau - 2} k_v^2 ) \tensor_{k_v, \tau} \Fpbar,$$
where $2 \leq k_\tau \leq p + 1$ and $0 \leq w_\tau \leq p - 1$, and the $w_\tau$ are not all $p - 1$.  Let $\Gamma = \prod_{v |p} \GL_2(k_v)$.  Then the irreducible $\Fpbar$-representations of $\Gamma$ are $\sigma = \tensor_{v |p} \sigma_v$ with $\sigma_v$ as above, and every weight factors through $\Gamma$.  We call the irreducible $\Fpbar$-representations of $\GL_2(k_v)$ {\emph {Serre weights at $v$}}.

Buzzard, Diamond, and Jarvis in \cite{BDJ} formulated a Serre-type conjecture for Hilbert modular forms in the case where $p$ is unramified in $F$.  We would like to have a conjecture in the general case.  
We may assume that $F \neq \Q$, as otherwise the conjecture is well-known (and mostly proved!).

Given a continuous, irreducible, totally odd Galois representation $\rho: \Gal(\overline{F}/F) \to \GL_2(\Fpbar)$, let $W(\rho)$ denote the set of weights for which it is modular; we explain below what is meant by ``modular.''  For each $v | p$ we will construct a set $W^?_v(\rho)$ of Serre weights at $v$ and conjecture that 
$$ W(\rho) = \left\{ \sigma = \bigotimes_{v | p} \sigma_v : \forall v,  \sigma_v \in W^?_v(\rho) \right\}.$$
This allows us to treat each $v | p$ separately.

In the next section we will state the conjecture in two equivalent forms, very much in the spirit of Florian Herzig's reformulation of the \cite{BDJ} conjecture.  The proof that they are equivalent (Theorems \ref{compthm} and \ref{compthmlvls}) relies heavily on Herzig's ideas in \cite{Florian}, \S 14.  In the third section we prove a theorem towards our conjecture; it shows, in many cases when the restriction of $\rho$ to a decomposition group at a place $v |p$ is irreducible, that the $v$-component of any modular weight does indeed lie in $W_v^?(\rho)$.  This statement, Theorem \ref{mainresult}, generalizes the main result of \cite{thesispreprint} and is proved by a similar argument; it was proved before the conjecture was formulated and played an important role in motivating it.  Finally, in the last section we use Demb\'{e}l\'{e}'s computations of Hilbert modular forms over $\Q(\sqrt{5})$ and their weights to obtain some computational evidence in support of the conjecture.

The author is grateful to Richard Taylor for encouraging him to study this question and for many enlightening conversations, to Lassina Demb\'{e}l\'{e} for making available the results of his computations of Hilbert modular forms over $\Q(\sqrt{5})$ and for performing new ones, and to Fred Diamond for useful conversations.

\section{A conjecture}
First we introduce the notion of modularity.  Let $D$ be a quaternion algebra over $F$ which is split at exactly one real place of $F$ and at all places over $p$.  Let $G = \Res_{F/\Q}(D^\ast)$ be the associated reductive group; for an open compact subgroup $U \subset G(\A^\infty)$ we have a Shimura curve $M_U /F$ whose complex points are
$$ M_U(\C) = G(\Q) \backslash G(\A^\infty) \times (\C - \R) / U.$$
The $M_U$ are not in general geometrically connected.  Let the abelian variety $\Pic^0(M_U) / F$ be the identity component of the relative Picard scheme of $M_U$, which parametrizes line bundles locally of degree zero.  

Let $U^\prime = \ker((D \tensor \hat{\Z})_p^\ast = \prod_{v |p} \GL_2(\mathcal{O}_v) \to \GL_2(\mathcal{O}_F /p))$, and let $U^{\prime \prime} = \ker( \prod_{v |p} \GL_2(\mathcal{O}_v) \to \prod_{v|p} \GL_2(k_v))$.  Clearly $U^\prime \subset U^{\prime \prime}$.  We say that an open compact $U \subset G(\A^\infty)$ is {\emph {of type $(\ast )$}} if $U = U^\prime \times U^p$, where $U^p \subset G(\A^{\infty, p})$.  Let $V = \prod_{v | p} \GL_2(\mathcal{O}_v) \times U^p$.  If $U^p$ is sufficiently small as in section 3.1 of \cite{thesispreprint}, then $M_U / M_V$ is a Galois cover with group $V / U = \GL_2(\mathcal{O}_F /p)$.  Hence we have an action of $V /U$ on $\Pic^0 (M_U)$.
\begin{definition}
An irreducible Galois representation $\rho: \Gal(\overline{F}/F) \to \GL_2(\Fpbar)$ is {\emph {modular of weight $\sigma$}} if there exists a quaternion algebra $D / F$ as above and an open compact $U \subset (D \tensor \hat{\Z})^\ast \subset G(\A^\infty)$ of type $(\ast )$, such that $(\Pic^0(M_U) [p] \tensor_{\Fpbar} \sigma)^{\GL_2(\mathcal{O}_F / p)} = (\Pic^0(M_{U^{\prime \prime} \times U^p}) [p] \tensor_{\Fpbar} \sigma)^{\Gamma}$ has $\rho$ as a Jordan-H\"{o}lder constituent. 
\end{definition}

Fix a place $\p | p$ of $F$; we will now define $W^?_\p(\rho)$.  Choose a decomposition subgroup $G_\p \subset \Gal(\overline{F}/F)$ at $\p$, and let $I_\p$ and $I_\p^\prime$ be the corresponding inertia and wild inertia subgroups.  Denote by $I_{t,\p} = I_\p / I_\p^\prime$ the tame inertia, and let the residue field $k_\p$ have cardinality $q = p^s$.

We will state our conjecture in the language of Herzig's reformulation of the \cite{BDJ} conjecture.  Let $I$ be the set of embeddings $k_\p \hookrightarrow \Fpbar$, and as in \cite{thesispreprint}, let $\tau_0, \dots, \tau_{s-1}$ be a labeling of its elements such that $\tau_{j - 1} = \tau_j^p$ for all $j \in \Z /s\Z$.  Similarly, let $k_\p^\prime$ be a quadratic extension of $k_\p$ and fix a labeling $\sigma_0, \dots, \sigma_{2s - 1}$ of the embeddings $k_\p^\prime \hookrightarrow \Fpbar$ such that $\sigma_{i-1} = \sigma_i^p$ for all $i \in \Z /2s\Z$ and such that $\sigma_i |_{k_\p} = \tau_{\pi(i)}$, where $\pi: \Z /2s\Z \to \Z /s\Z$ is the natural projection.  Given such an embedding $\tau \in I$ (resp. $\sigma: k_\p^\prime \hookrightarrow \Fpbar$), let $\lambda_\tau: I_{t,\p} \simeq \varprojlim \F_{p^n}^\ast \to \Fpbar$ (resp. $\psi_\sigma: I_{t,\p} \to \Fpbar$) be the corresponding fundamental character of level $s$ (resp. $2s$).  Often we write $\lambda_j, \psi_i$ for $\lambda_{\tau_j}, \psi_{\sigma_i}$.  Note that Herzig's convention is $\psi_{i+1} = \psi_i^p$; the reader should bear this in mind when comparing our work with his.

If $b = \sum_{j = 0}^{s-1} w_j p^{s-j}$ and $a - b = \sum_{j = 0}^{s-1} (k_j - 2) p^{s-j}$ for $0 \leq w_j \leq p - 1$ and $2 \leq k_j \leq p + 1$, then we denote
$$ F(a,b) = \bigotimes_{j \in \Z /s\Z} (\det\nolimits^{w_j} \Sym^{k_j - 2} k_\p^2) \tensor_{k_\p, \tau_j} \Fpbar .$$
Of course this notation comes from the theory of Weyl modules, but for the purposes of this article we may take the expression above as a definition.

Given $\rho |_{I_\p}$, we first associate to it a characteristic zero representation of $\GL_2(k_\p)$ as in \cite{Florian}, Def. 14.1.  Here $I(\chi_1, \chi_2)$ are the usual principal series, while the $\Theta(\xi)$, for $\xi: k_\p^\prime \hookrightarrow \Fpbar$, are the cuspidal representations (see, for instance, \cite{DeligneLusztig}).
\begin{definition} 
\begin{enumerate}
\item If $\rho |_{I_\p} \sim \left( \begin{array}{cc} \prod_{j \in \Z /s\Z} \lambda_j^{m_j} & 0 \\ 0 & \prod_j \lambda_j^{n_j} \end{array} \right)$, then $V_\p(\rho) = I(\prod \tau_j^{m_j}, \prod \tau_j^{n_j}) $.
\item If $\rho |_{I_\p} \sim \left( \begin{array}{cc} \prod_{i \in \Z /2s\Z} \psi_i^{m_i} & 0 \\ 0 & \prod_i \psi_i^{m_{i+s}} \end{array} \right) \prod_{j \in \Z /s\Z} \lambda_j^{w_j}$, then $V_\p(\rho) = \Theta(\prod \sigma_i^{m_i}) \tensor \prod_j \tau_j^{w_j}$.
\end{enumerate}
\end{definition}

Since $V_\p(\rho)$ can be realized over $\overline{\Z}_p$, we may consider its reduction modulo $p$, denoted $\overline{V_\p(\rho)}$.  For any representation $V$, we write $JH(V)$ for the set of its Jordan-H\"{o}lder constituents.  The sets $JH(\overline{V_\p(\rho)})$ are computed in \cite{Fredpreprint}.

In Lemma \ref{centchar} we compute the determinant of $\rho |_{I_\p}$, and hence the central character of any modular weight.  If $e \geq p$, we conjecture that {\emph {all}} weights with this central character are modular.  Indeed, this is suggested by the fact that we already conjecture this ``maximal'' set of weights when $e = p - 1$, as can be seen from Theorems \ref{compthm} and \ref{compthmlvls}, and by the observation that the number of conjectured modular weights increases with $e$ for $e \leq p - 1$.

Let $Y_\p$ be the set of Serre weights at $\p$.  If $e \leq p -1$, let $\delta \in \Delta = [0, e - 1]^I$ be a vector whose components are choices of an integer $0 \leq \delta_\tau \leq e - 1$ for each $\tau \in I$.  Given $\delta$, we will define a multi-valued function $\mathcal{R}_\p^\delta : Y_\p \to Y_\p$ for which we conjecture the following:

\begin{conjecture} \label{mainconj}
Let $\rho: \Gal(\overline{F} /F) \to \GL_2(\Fpbar)$ be continuous, irreducible, totally odd, and tame at $\p$.  Then
\begin{enumerate}
\item $ W_\p^?(\rho) = \bigcup_{\delta \in \Delta} \mathcal{R}_\p^\delta (JH(\overline{V_\p(\rho)}))$ if $e \leq p -1$.
\item $ W_\p^?(\rho) = \left\{ F(a,b): \det \rho |_{I_\p} = \lambda_0^{a + b + \sum_{j = 0}^{s-1} e p^j} \right\}$ if $e \geq p$.
\end{enumerate}
\end{conjecture}

We will now assume $e \leq p - 1$, fix $\delta \in \Delta$, and construct the map $\mathcal{R}_{\p}^{\delta}$.  Given $F(a,b)$, define $\alpha(j) = p + 1 - k_j \in [0, p-1]$ for every $j \in \Z /s\Z$.  Define $x_j$ to be the integer such that $\alpha(j) + x_j p \in [1 + 2 \delta_j - (e-1), p + 2 \delta_j - (e - 1)]$.  Under the assumption that $e \leq p -1$, we have $x_j \in \{ -1, 0, 1 \}$ for all $j$.  We say that $F(a,b)$ is a $\delta$-{\emph{regular}} Serre weight at $\p$ if the $x_j$ are all zero.  If $F(a,b)$ is not $\delta$-regular, then for every $j \in \Z /s\Z$ we define $\theta_j = x_{j + n}$, where $n$ is the smallest positive integer such that $x_{j + n} \neq 0$.  

Suppose first that $F(a,b)$ is a $\delta$-regular Serre weight.  Then we define
$$ \mathcal{R}_\p^\delta(F(a,b)) = \left\{ F(c,d) : {{c \equiv b - \sum_{j = 0}^{s-1} (1 + \delta_j) p^{s-j} \mod p^s - 1} \atop  {d \equiv a - \sum_{j = 0}^{s-1} (e-1 - \delta_j) p^{s-j} \mod p^s - 1}} \right\}.$$

If $F(a,b)$ is irregular, things become more complicated.  We define a collection $\mathcal{S}^\delta(F(a,b))$ of subsets of $ \Z / s\Z$ as follows.  Let $S \subset \Z /s\Z$.  Then $S \in \mathcal{S}^\delta(F(a,b))$ if and only if for every $j \in S$ the following two conditions hold:
\begin{enumerate}
\item  One of the following two conditions holds:
\begin{enumerate}
\item $x_j = -1$ or $\alpha(j) \in [2 \delta_j - (e-1), p - 1 + 2 \delta_j - (e-1)] \cap [0,p-1]$, and there is an integer $n \geq 0$ such that $x_{j + m} = 1 + 2 \delta_{j+m} - (e - 1)$ for $1 \leq m \leq n$ (if any such $m$ exists) and $x_{j+n+1} = 1$.
\item $x_j = 1$ or $\alpha(j) \in [2 + 2 \delta_j - (e-1), p + 1 + 2 \delta_j - (e-1)] \cap [0, p - 1]$, and there is an $n \geq 0$ such that $x_{j + m} = p + 2 \delta_{j+m} - (e-1)$ for all $1 \leq m \leq n$ and $x_{j + n + 1} = -1$.
\end{enumerate}
\item In either of the cases above, $j + m \not\in S$ for $1 \leq m \leq n$.
\end{enumerate}

We emphasize that $\mathcal{S}^\delta(F(a,b))$ depends only on $a - b$.  Finally, writing $F = F(a,b)$, we can give the general definition:
$$ \mathcal{R}_\p^\delta(F) = \left\{ F(c,d) : {{c \equiv b + \sum_{j \in S} \theta_j p^{s-j} - \sum_{j = 0}^{s-1} (1 + \delta_j) p^{s-j} \mod p^s - 1} \atop  {d \equiv a - \sum_{j = 0}^{s-1} (e-1 - \delta_j) p^{s-j} - \sum_{j \in S} \theta_j p^{s-j} \mod p^s - 1}} : S \in \mathcal{S}^\delta(F) \right\}.$$

\begin{lemma} [\cite{Florian}, Lemma 14.3] \label{florianlemma}
If $\rho$ is of level $2s$, then 
$ \sigma_\p = \bigotimes_{\tau \in I} (\det\nolimits^{w_\tau} \Sym^{k_\tau - 2} k_\p^2 ) \tensor_{k_\p, \tau} \Fpbar $
is a Jordan-H\"{o}lder constituent of $\overline{V_\p(\rho)}$ if and only if for each $\tau \in I$ there is a labeling $\{ \tilde{\tau}, \tilde{\tau}^\prime \}$ of its two lifts to $k_\p^\prime$ such that 
$$ \rho |_{I_\p} \sim \prod_{\tau} \lambda_\tau^{w_\tau + k_\tau - 2} \left( \begin{array}{cc} \prod_\tau \psi_{\tilde{\tau}}^{p + 1 - k_\tau} & 0 \\ 0 & \prod_\tau \psi_{\tilde{\tau}^\prime}^{p + 1 - k_\tau} \end{array} \right) .$$
\end{lemma}

\begin{theorem} \label{compthm}
Suppose that $\rho |_{I_\p}$ is of level $2s$.  Then $W_\p^?(\rho)$ consists precisely of those Serre weights at $\p$
\begin{equation} \label{wtp}
 \sigma_\p = \bigotimes_{\tau \in I} (\det\nolimits^{w_\tau} \Sym^{k_\tau - 2} k_\p^2 ) \tensor_{k_\p, \tau} \Fpbar 
 \end{equation}
such that for each $\tau \in I$ there exists a labeling $\{ \tilde{\tau}, \tilde{\tau}^\prime \}$ of its two lifts to $k_\p^\prime$ and an integer $0 \leq \delta_\tau \leq e - 1$ such that 
$$ \rho |_{I_\p} \sim \prod_{\tau \in I} \lambda_\tau^{w_\tau} \left( \begin{array}{cc} \prod_\tau \psi_{\tilde{\tau}}^{k_\tau - 1 + \delta_\tau} \psi_{\tilde{\tau}^\prime}^{e - 1 - \delta_\tau} & 0 \\
0 & \prod_\tau \psi_{\tilde{\tau}}^{e - 1 - \delta_\tau} \psi_{\tilde{\tau}^\prime}^{k_\tau - 1 + \delta_\tau} \end{array} \right) .$$
\end{theorem}
\begin{proof}
If $e \geq p$ the theorem is evident, so we assume from now on that $e \leq p - 1$.
Let $L_\p^\delta(\rho)$ be the set of weights satisfying the condition in the statement above for a given $\delta$.  We claim that $\mathcal{R}_\p^\delta(JH(\overline{V_\p(\rho)})) = L_\p^\delta(\rho)$ for every choice of $\delta$.

Fix $\delta$, and suppose that $F = F(a,b)$ is a Jordan-H\"{o}lder constituent of $\overline{V_\p(\rho)}$.  Without loss of generality, we may assume that $b = 0$, and we write
$a = \sum_{j = 0}^{s-1} a_j p^j$ with $0 \leq a_j \leq p - 1$.  Set $\alpha(j) = p - 1 -a_j$.  Then by Lemma \ref{florianlemma} we have
\begin{eqnarray} \label{matrixexp}
 \rho |_{I_\p} & \sim &\prod_{j \in \Z /s\Z}  \lambda_j^{a_j} \left( \begin{array}{cc} \prod_{i \in J} \psi_i^{\alpha(i)} & 0 \\ 0 & \prod_{i \in J^c} \psi_i^{\alpha(i)} \end{array} \right) = \\ \nonumber
& & \prod_{j \in \Z/s\Z} \lambda_j^{a_j + p - e + \delta_j} \left( \begin{array}{cc} \prod_{i \in J} \psi_i^{\alpha(i) + e - 1 - \delta_i} \psi_{i + s}^{e - 1 - \delta_i} & 0 \\ 0 & \prod_{i \in J^c} \psi_i^{\alpha(i) + e - 1 - \delta_j} \psi_{i + s}^{e - 1 - \delta_j} \end{array} \right) ,
\end{eqnarray}
where $J \subset \Z /2s\Z$ is a subset such that $|J| = s$ and $\pi(J) = \Z /s\Z$ and we write $\alpha(i)$ for $\alpha(\pi(i))$ and $\delta_i$ for $\delta_{\pi(i)}$.  Our goal now is to write $\prod_{i \in J} \psi_i^{\alpha(i)}$ in the form $\eta \prod_{i \in J} \psi_i^{\beta(i)}$, where $\beta(i) \in [ 1 + 2 \delta_i - (e - 1), p + 2 \delta_i - (e - 1)]$ and $\eta$ is a character of level $s$; from such an expression we will read off a weight in $L_\p^\delta(\rho)$.

Observe first that such an expression is unique if it exists.  Indeed, suppose that $\eta \prod_J \psi_i^{\beta(i)}$ and $\eta^\prime \prod_J \psi_i^{\beta(i)^\prime}$ are two such expressions.  Then $\psi = \prod_{i \in J} \psi_i^{\beta(i) - \beta(i)^\prime}$ is a character of level $s$, whence
$$ \psi^{1 - p^s} = \prod_{i \in J} \psi_i^{\beta(i) - \beta(i)^\prime} \psi_{i + s}^{\beta(i)^\prime - \beta(i)} = 1.$$
Since $| \beta(i) - \beta(i)^\prime | \leq p - 1$ for all $i \in J$, it is evident that we must have $\beta(i) = \beta(i)^\prime$ for all $i$.

Two issues must be dealt with in obtaining the desired expression.  First, $J$ is not specified uniquely by $(\ref{matrixexp})$.  Indeed, if $\alpha(j) = 0$ for some $j \in \Z /s\Z$, then we can choose either of the elements of $\pi^{-1}(j)$ to lie in $J$.  In this case, $\overline{V_\p(\rho)}$ has fewer constituents than usual, but the ambiguity in $J$ allows us to produce several modular weights from each constituent.  Second, the $\alpha(j)$ need not lie in the range $[1 + 2 \delta_j - (e-1), p + 2 \delta_j - (e-1)]$, so we must  ``carry''  exponents.  If $e = 1$, this problem occurs only when $\alpha(j) = 0$, which is exactly when the first problem arises.  In general we do not have this coincidence, whence the relative complexity of our construction to the analogous one in \cite{Florian}, \S 14.  If $e = 1$, then the argument below reduces precisely to Herzig's argument.

For each $j \in \Z /s\Z$, let $x_j \in \Z$ be such that $\alpha(j) + x_j p \in [ 1 + 2 \delta_j - (e - 1), p + 2 \delta_j - (e - 1)]$.  For instance, if $e = 1$, then $x_j = 1$ if $\alpha(j) = 0$ and $x_j = 0$ otherwise.  Observe that if $e \leq p - 1$, then $x_j \in \{ -1, 0, 1 \}$; moreover, $\alpha(j) + x_j (p -1)$ also lies in the specified range.  

As in \cite{Florian}, we consider an interval in $\Z /s\Z$ to be a sequence $[[ j, n ]] = \{ j, j+1, \dots, n \}$.  The predecessor of $[[j,n]]$ is $j - 1$ (these correspond, of course, to Herzig's successors; the difference is a consequence of our opposite conventions).  The terminus of $[[j,n]]$ is $n$.  We define $\mathcal{L}_\delta$ as the set of all pairs $(\alpha, \mathcal{I})$, where $\alpha: \Z /s\Z \to [ 0, p - 1]$ is a map, $\mathcal{I}$ is a collection of disjoint intervals, each labeled with a sign, and the following axioms are satisfied.  For each $j \in \Z /s\Z$, given $\alpha$, we can formally define $x_j$ as above.  If $j - 1$ is the predecessor of an interval and $n$ is its terminus, and $x_n \neq 0$, then define $z_{j - 1} = x_n$.  Otherwise it will follow from the third axiom that $n + 1 \in \bigcup \mathcal{I}$ and we define $z_{j - 1} = z_n$.  Thus $z_j = \pm 1$.  Finally, let $C^\pm = \{ j \in \Z /s\Z : x_j = \pm 1 \}$.  The axioms are:
\begin{enumerate}
\item For each interval $I \in \mathcal{I}$, either $I \subset C^+ \cup \{ j : \alpha(j) = 1 + 2\delta_j - (e - 1)\}$ or $I \subset C^- \cup \{ j : \alpha(j) = p + 2\delta_j - (e - 1) \}$.
\item If $j \in \bigcup \mathcal{I}$ and $\alpha(j) \neq 0$, then $j$ is the terminus of its interval.
\item If $j \in \bigcup \mathcal{I}$, then $\alpha(j) \in \{ 1 + 2\delta_j - (e - 1), p + 2 \delta_j - (e - 1) \}$ if and only if $j$ is the terminus of an $\mathcal{I}$-interval and the predecessor of a negative $\mathcal{I}$-interval.  
\item If $j \not\in \bigcup \mathcal{I}$ and $x_j \neq 0$, then $ j + 1 \in \bigcup \mathcal{I}$.
\item If a positive $\mathcal{I}$-interval has predecessor $j$, then either $j$ lies in an $\mathcal{I}$-interval and satisfies $z_j \neq x_j$, or $j$ does not lie in any interval and $\alpha(j) \in [1 - z_j + 2 \delta_j - (e - 1), p - z_j + 2 \delta_j - (e-1)]$.
\item If a negative $\mathcal{I}$-interval $I$ has predecessor $j$, then either $j$ lies in an $\mathcal{I}$-interval and satisfies $z_j = x_j$, or $j$ lies in an $\mathcal{I}$-interval and $\alpha(j) = 1 + 2 \delta_j - (e-1)$ (resp. $p + 2 \delta_j - (e - 1))$ if $z_j = 1$ (resp. $z_j = -1$), or else $j$ does not lie in any interval and $\alpha(j) \in [1 + z_j + 2 \delta_j - (e-1), p + z_j + 2 \delta_j - (e-1)]$.
\end{enumerate}

Similarly, let $\mathcal{M}_\delta$ be the set of pairs $(\beta, \mathcal{I})$, where $\mathcal{I}$ as before is a collection of signed intervals and $\beta: \Z /s\Z \to \Z$ is a map such that for every $j \in \Z /s\Z$, we have $\beta(j) \in [1 + 2 \delta_j - (e - 1), p + 2 \delta_j - (e - 1)]$.  Let $y_j$ be the integer such that $\beta(j) - y_j p \in [0, p - 1]$.  Let $D^\pm = \{ j \in \Z /s\Z: y_j = \pm 1 \}$.  If $j$ is the predecessor of an interval, let $u_j$ be the number defined in the same way as $z_j$, but with $x_j$ replaced by $y_j$ in the definition.  We require that $(\beta, \mathcal{I})$ satisfy the following axioms:
\begin{enumerate}
\item For each interval $I \in \mathcal{I}$, either $I \subset D^+ \cup \{j: \beta(j) = p - 1 \}$ or $I \subset D^- $.
\item The set of termini of $\mathcal{I}$-intervals is $D^+ \cup D^-$. 
\item If a positive $\mathcal{I}$-interval has predecessor $j$, then either $j$ lies in an $\mathcal{I}$-interval and satisfies $u_j \neq y_j$, or $j$ does not lie in any interval and $\beta(j) \in [u_j, p - 1 + u_j]$.
\item If a negative $\mathcal{I}$-interval has predecessor $j$, then either $j$ lies in an $\mathcal{I}$-interval and satisfies $u_j = y_j$, or $j$ does not lie in any interval and $\beta(j) \in [-u_j, p - 1 - u_j]$.
\end{enumerate}

There is a bijection $\xi: \mathcal{L}_\delta \to \mathcal{M}_\delta$ which can be written down as follows.  Like Herzig, we represent the function $\alpha$ by the string of numbers $\alpha(0), \alpha(1), \dots, \alpha(s-1)$.  We underline each $\mathcal{I}$-interval and put its sign after its last entry.  Pairs $(\beta, \mathcal{I})$ are written similarly.  Then $\xi$ acts as follows, where $j$ is always the predecessor of the last interval, in the third line $k$ is the predecessor of the first interval, and we assume $x_j = 0$ in the second line and $x_j \neq 0, x_k \neq 0$ in the third:
\begin{eqnarray*}
x^\prime, \underline{(0, \dots, 0,) x}_\pm & \mapsto & x^\prime \pm z_j, \underline{(z_j(p-1), \dots, z_j (p-1),) x + z_j p}_\pm \\
y^\prime, \underline{(0 \dots 0,) y,}_\pm \underline{(0 \dots 0,) y^{\prime \prime}}_- & \mapsto &
y^\prime \pm z_j, \underline{(z_j(p-1) \dots z_j(p-1)), y + z_j(p-1),}_\pm \underline{(z_j(p-1) \dots, }_- \\
w^\prime, \underline{(0 \dots 0) w,}_\pm \underline{(0 \dots 0) w^{\prime \prime}}_{\pm^\prime} & \mapsto &
w^\prime \pm z_k, \underline{(z_k(p-1), \dots ), w \pm z_k p \pm^\prime z_j}_\pm \underline{(z_j(p-1), \dots), w^{\prime \prime} + z_j p}_{\pm^\prime}
\end{eqnarray*}

All other entries are unchanged by $\xi$.  The reader may verify that $\xi$ is indeed a bijection between $\mathcal{L}_\delta$ and $\mathcal{M}_\delta$.  It does not affect the collection $\mathcal{I}$ of signed intervals.  We will prove below (Lemma \ref{pospred}) that if $S \subset \Z /s\Z$, then $S \in \mathcal{S}^\delta(F(a, 0))$ if and only if $S$ is the set of predecessors of positive intervals in $\mathcal{I}$ for some $(\alpha, \mathcal{I}) \in \mathcal{L}_\delta$ (hence for some $(\beta, \mathcal{I}) \in \mathcal{M}_\delta)$, where $\alpha$ is derived from $a$ as above.

Given $S \in \mathcal{S}^\delta(F)$, let $(\alpha, \mathcal{I}) \in \mathcal{L}_\delta$ be such that $S$ is the set of predecessors of positive intervals.  Let $J_+$ (resp. $J_-$) be the elements of $J$ whose projections to $\Z /s\Z$ are predecessors of positive (resp. negative) intervals, and similarly for $J^c$.  Let $\tilde{\mathcal{I}}$ be the collection of intervals in $\Z /2s\Z$ that project to $\mathcal{I}$-intervals, and let $J_0$ be the elements of $J$ that do not lie in any $\tilde{\mathcal{I}}$-intervals.  Then as in \cite{Florian} we observe that $\prod_{i \in J} \psi_i^{\alpha(i)} = \chi \prod_{i \in J_+ \cup J_+^c} \psi_i^{- z_i}$, where
$$
\chi = \prod_{i \in J_+} \psi_i^{\alpha(i) + z_i} \prod_{J_0 \backslash (J_+ \cup J_-)} \psi_i^{\alpha(i)} \prod_{J_-} \psi_i^{\alpha(i) - z_i} \prod_{i - 1 \in J_- \cap J^c_+ \atop [[i,n]] \in \tilde{\mathcal{I}} } (\psi_i^{p-1} \psi_{i+1}^{p-1} \cdots \psi_n^{p})^{z_{i-1}} \prod_{J \backslash (J_+ \cup J_- \cup J_0)} \psi_i^{\alpha(i)}$$
and it is not hard to see that in this expression, one of every pair $\{ \psi_i, \psi_{i + s} \}$ appears with exponent zero and the other appears with exponent in the range $[ 1 + 2 \delta_j - (e-1), p + 2 \delta_j - (e-1)]$.  Hence 
$$ \rho |_{I_\p} \sim \prod_{j \in \Z /s\Z} \lambda_j^{a_j + p - e + \delta_j} \left( \begin{array}{cc} \chi_1 & 0 \\ 0 & \chi_2 \end{array} \right) \prod_{j \in S} \lambda_j^{-z_j} , $$
where each $\psi_i$, $i \in \Z /2s\Z$, appears with exponent $e -1 - \delta_i$ in one of $\chi_1, \chi_2$ and with some exponent $\beta(i) = \beta_{\pi(i)}$ in the range $[1 + \delta_j, p + \delta_j]$ in the other.  From such an expression we can read off a weight $F(A,B) \in L_\p(\rho)$.

Now, from $(\ref{matrixexp})$ we see that $\det \rho |_{I_\p} = \prod_{j \in \Z/s\Z} \lambda_j^{a_j} = \lambda_0^{\sum_{m = 0}^{s-1} a_{s-j} p^j}$.  Let $1_S: I \to \{ 0, 1 \}$ be the characteristic function of $S$.  Then from the displayed expressions above we find that
$$ \det \rho |_{I_\p} = \prod_{j \in \Z /s\Z} \lambda_j^{2a_j - (e-1) + \delta_j + \beta_j} \prod_{j \in S} \lambda_j^{-2z_j} = \lambda_0^{\sum_{m=0}^{s-1} (2a_{-m} - (e-1) + \delta_{-m} + \beta_{-m} - 2 \cdot 1_S(-m) z_{-m}) p^j}.
$$
Hence, noting that for $j \in S$ we have $w_j = z_j$, we find that
\begin{eqnarray*}
B & \equiv & a + \sum_{m=0}^{s-1} (\delta_{j} - (e-1)) p^{s-j} - \sum_{j \in S} w_j p^{s-j} \mod p^s - 1 \\
A & \equiv & \sum_{j \in S} w_j p^{s-j} - \sum_{j = 0}^{s-1} (\delta_j + 1) p^{s-j} \mod p^s - 1
\end{eqnarray*}

It remains to check that any other weight $F(\tilde{A}, \tilde{B})$ satisfying the same congruences is also contained in $L_\p^\delta(\rho)$.  The only cases when more than one weight satisfies such a congruence are the pairs $F(b,b), F(p^s - 1 + b, b)$ for some $b$.  But then it is obvious from the definition of $L_\p^\delta(\rho)$ that one of these weights is contained there if and only if the other one is.  Hence we have shown that $\mathcal{R}_\p^\delta(F) \subset L_\p^\delta(\rho)$.

Conversely, suppose that $F(a,b) \in L_\p^\delta(\rho)$.  We may assume without loss of generality that $b = 0$, and as usual write $a = \sum_{j=0}^{s-1} a_{j} p^{s-j}$.  Then,
$$ \rho |_{I_\p} \sim \left( \begin{array}{cc} \prod_{i \in L} \psi_i^{\beta(i)} & 0 \\ 0 & \prod_{i \in L^c} \psi_i^{\beta(i)} \end{array} \right) \prod_{i \in \Z /2s\Z} \psi_i^{e-1-\delta_i} ,$$
where $L \subset \Z / 2s\Z$ is mapped bijectively to $\Z /s\Z$ by $\pi$ and $ \beta(i) =  a_{\pi(i)} + 1 + 2 \delta_i - (e-1) \in [ 1 + 2\delta_i - (e - 1), p + 2 \delta_i - (e - 1)]$.  Let $y_i$ be an integer such that $\beta(i) - y_i p \in [0, p-1]$; under our assumptions on $e$, we have $y_i \in \{ -1, 0, 1\}$.  Let $D^\pm = \{ i \in \Z /2s\Z : y_i = \pm 1 \}$.  We now define a collection $\mathcal{I}$ of intervals in bijection with $D^+ \cup D^-$ as follows.  If $i \in D^+$ and $i \in L$ (resp. $i \in L^c$), choose $n$ such that $[[ n, i ]] \subset L$ (resp. $L^c$) and $\beta(m) = p - 1 $ for all $m \in [[ n, i ]] \backslash \{ i \}$, and such that $n$ is minimal for this property (i.e. $n - 1$ will not work).  Then $[[ n, i ]]$ is the interval corresponding to $i$, and we let it be negative if and only if $y_{n - 1} = 1$ or $y_{n-1} = 0$ and $n - 1 \in L$ (resp. $L^c$).

Similarly, if $i \in D^-$, then the corresponding interval is $[[ i ]]$.  It is negative if and only if $y_{i - 1} = -1$ or $y_{i - 1} \in \{ 1 + 2 \delta_{i-1} - (e-1), p + 2 \delta_{i-1} - (e-1) , p - 1 \}$.

It is easy to see that $(\beta, \mathcal{I}) \in \mathcal{M}_\delta$.  Let $L_+$, $L_-$, $L_0$ be defined as before, and let $S = L_+ \cup L_+^c$ be the set of predecessors of positive intervals.  We invert the previous construction to find that 
$$ \rho |_{I_\p} \sim \left( \begin{array}{cc} \chi_1 & 0 \\ 0 & \chi_2 \end{array} \right) \prod_{j \in \Z /s\Z} \lambda_j^{e-1-\delta_j} \prod_{j \in S} \lambda_j^{u_j}, $$
where
$$ \chi_1 = \prod_{i \in L_+} \psi_i^{\beta(i) - u_i} \prod_{L_0 \backslash (L_+ \cup L_-)}^{\beta(i)} \prod_{L_-} \psi_i^{\beta(i) + u_i} \prod_{i - 1 \in L_- \cup L_+^c \atop [[ i, n]] \in \tilde{\mathcal{I}}} (\psi_i^{p-1} \cdots \psi_n^p)^{-u_{i-1}} \prod_{L \backslash (L_+ \cup L_- \cup L_0)} \psi_i^{\beta(i)}, $$
and $\chi_2$ is the same but with the roles of $L$ and $L^c$ reversed.  Each $\psi_i$ appears with non-zero exponent in at most one of $\chi_1, \chi_2$, and this exponent always lies in the range $[1, p-1]$.  Thus we have obtained an expression of the form
$$ \rho |_{I_\p} \sim \left( \begin{array}{cc} \prod_{i \in L} \psi_i^{\alpha(i)} & 0 \\ 0 & \prod_{i \in L^c} \psi_i^{\alpha(i)} \end{array} \right) \prod_{j \in \Z /s\Z} \lambda_j^{e-1-\delta_j} \prod_{j \in S} \lambda_j^{u_j}.$$

Using Lemma \ref{florianlemma} we can read off a weight $F(A,B) \in JH(\overline{V_\p(\rho)})$.
Moreover, clearly $(\alpha, \mathcal{I}) = \xi^{-1} (\beta, \mathcal{I})$, whence $S \in \mathcal{S}^\delta(F(A,B))$.  Comparing two expressions for $\det \rho |_{I_\p}$ as before, we find that
$$ \sum_{j = 0}^{s-1} (a_j + e) p^{s-j} \equiv \sum_{j = 0}^{s-1} (\alpha(j) + 2[(e -1 - \delta_j) + 1_S(j) u_j ])p^{s-j} \mod p^s -1.$$
Clearly $u_j = w_j$ for $j \in S$.  Also we see that
\begin{eqnarray*}
B & \equiv & \sum_{j = 0}^{s-1} (e-1 - \delta_j + 1_S(j) u_j + \alpha(j))p^{s-j} \equiv \sum_{j=0}^{s-1} (a_j + 1 + \delta_j) p^{s-j}  - \sum_{j \in S} u_j p^{s-j} \mod p^s -1 \\
A & \equiv & \sum_{j = 0}^{s-1} (e-1 - \delta_j + 1_S(j) u_j) p^{s-j} \mod p^s -1
\end{eqnarray*}

Hence, $F(a,b) \in \mathcal{R}_\p^\delta(F(A,B))$.  This completes the proof that $\mathcal{R}_\p^\delta(JH(\overline{V_\p(\rho)})) = L_\p^\delta (\rho)$.
\end{proof}
A very similar argument establishes an analogous statement in the level $s$ case:
\begin{theorem} \label{compthmlvls}
Suppose that $\rho |_{I_\p}$ is of level $s$ and, as always, tame at $\p$.  Then $W_\p^?(\rho)$ consists precisely of the Serre weights at $\p$ as in $(\ref{wtp})$ for which there exist a set $J \subset I$ and an integer $0 \leq \delta_\tau \leq e - 1$ for each $\tau \in I$ such that 
$$ \rho |_{I_\p} \sim \prod_{\tau \in I} \lambda_\tau^{w_\tau} \left( \begin{array}{cc} \prod_{\tau \in J} \lambda_\tau^{k_\tau - 1 + \delta_\tau} \prod_{\tau \not\in J} \lambda_\tau^{e - 1 - \delta_\tau} & 0 \\
0 & \prod_{\tau \in J} \lambda_\tau^{e -1 - \delta_\tau} \prod_{\tau \not\in J} \lambda_\tau^{k_\tau - 1 + \delta_\tau} \end{array} \right) .$$
\end{theorem}

Finally we establish a lemma that was needed in the proof of Theorem \ref{compthm}.  
\begin{lemma} \label{pospred}
Let $\alpha : \Z /s\Z \to [0, p-1]$ be a function, and let $S \subset \Z /s\Z$.  Then $S \in \mathcal{S}(F(a,b))$ for some (hence all) weights $F(a,b)$ such that $a - b = \sum_{j = 0}^{s-1} (p - 1 - \alpha(j)) p^{s-j}$ if and only if $S$ is the set of predecessors of positive $\mathcal{I}$-intervals for some $(\alpha, \mathcal{I}) \in \mathcal{L}_\delta$.
\end{lemma}
\begin{proof}
It is easy to see from the axioms of $\mathcal{L}_\delta$ that the set of predecessors of positive intervals of any $(\alpha, \mathcal{I})$ lies in $\mathcal{S}^\delta(F(a,b))$.  

Conversely, suppose $S \in \mathcal{S}^\delta(F(a,b))$; we will construct an appropriate $\mathcal{I}$.  We let $j \in \bigcup \mathcal{I}$ if and only if there exists $n \geq 0$ such that $x_{j+n+1} = 0$ and for all $1 \leq m \leq n$ we have $j + m \not\in S$ and either $x_{j + m} = 1 + 2 \delta_{j + m} - (e - 1)$ for all $m$ or $x_{j + m} = p + 2 \delta_{j + m} - (e - 1)$ for all $m$.  We let $j \in \bigcup \mathcal{I}$ be the terminus of an $\mathcal{I}$-interval if and only if $j + 1 \not\in \bigcup \mathcal{I}$, or if $j + 1 \in \bigcup \mathcal{I}$ and $\alpha(j) \neq 0$, or $\alpha(j) = 0 \in \{ 1 + 2 \delta_j - (e -1), p + 2 \delta_j - (e - 1) \}$.  This specifies $\mathcal{I}$, and we define an $\mathcal{I}$-interval to be positive if and only if its predecessor is contained in $S$.  The reader may verify that $(\alpha, \mathcal{I}) \in \mathcal{L}_\delta$.
\end{proof}

\section{A theorem towards the conjecture}
As before, let $I$ be the set of embeddings $\tau: k_\p \hookrightarrow \Fpbar$.
Suppose the Galois representation $\rho: \Gal(\overline{F}/F) \to \GL_2(\Fpbar)$ is modular of a weight $\sigma$ whose $\p$-component is 
\begin{equation} \label{sigmap}
\sigma_\p = \bigotimes_{\tau \in I} (\det\nolimits^{w_\tau} \Sym^{k_\tau - 2} k_\p^2 ) \tensor_{k_\p, \tau} \Fpbar .
\end{equation}
Suppose that the restriction of $\rho$ to the decomposition subgroup $G_\p$ is irreducible.  Then as in \cite{thesispreprint} we have 
$$ \rho |_{I_\p}^{ss} \sim \left( \begin{array}{cc} \phi & 0 \\ 0 & \phi^q \end{array} \right),$$
where $\phi: I_{t,\p} = I_\p / I_\p^\prime \to \Fpbar^\ast$ is a character of level $2s$.  
Let $K$ be the maximal unramified extension of $F_\p$, and let $K^\prime /K$ be the totally ramified extension such that $\Gal(K^\prime /K) \simeq k_\p^\ast$.

The present argument is very similar to the one in \cite{thesispreprint}, so we refer the reader to that article and only indicate the differences.  In particular, the first four sections of \cite{thesispreprint} do not depend on the assumption that $p$ is unramified in $F$, so they hold in our case as well.  Suppose that $\rho$ is modular of weight $\sigma = \sigma_\p \tensor (\tensor_{v \neq \p} \sigma_v)$ and that $\sigma_\p$ is a Jordan-H\"{o}lder constituent of $Ind_{B}^{\GL_2(k_\p)} \theta$, where $B \subset \GL_2(k_\p)$ is the subgroup of upper triangular matrices and $\theta: B \to \Fpbar^\ast$ is given by
\begin{equation} \label{deftheta}
\theta: \left( \begin{array}{cc} a & b \\ 0 & d \end{array} \right)  \mapsto  \prod_{\tau: k_\p \hookrightarrow \Fpbar} \tau(ad)^{w_\tau} \tau(d)^{k_\tau - 2}.
\end{equation}

\begin{lemma} \label{centchar}
Write $k_j$ for $k_{\tau_j}$.  Then,
$$ \phi^{q + 1} = \prod_{i \in \Z /2s\Z} \psi_i^{2w_{\pi(i)} + k_{\pi(i)} - 2 + e}.$$
\end{lemma}
\begin{proof}
By \cite{thesispreprint}, Prop.~ 3.19, for all $\sigma \in \Gal(\overline{F}/F)$ we have $\det \rho(\sigma) = \chi(\sigma) \langle \sigma \rangle^{-1}$, where $\chi$ is the mod $p$ cyclotomic character and $\langle \cdot \rangle$ is the diamond operator map.  If $\sigma \in \Gal(\overline{K}/K) = I_\p$, suppose its image in $\Gal(K^\prime /K)$ is sent by the Artin reciprocity map to $j(\sigma) \in \mathcal{O}_\p^\ast / (1 + \p)$.  Then we have
$$ \phi^{q + 1}(\sigma) = \det \rho(\sigma) = \chi(\sigma) \langle \sigma \rangle^{-1} = \prod_{\tau: k_\p \hookrightarrow \Fpbar} \tau(j(\sigma))^{k_\tau - 2} \tau(j(\sigma))^{e} = \prod_{i \in \Z /2s\Z} \psi_i(\sigma)^{k_{\pi(i)} - 2 + e},$$
just as in the proof of \cite{thesispreprint}, Lemma 5.1.
\end{proof}

Assume from now on that $e \leq p - 1$.  Let $\boldsymbol{\mu} \in \Z^{s}$ be the vector whose components are given by $\mu_i = a_i + a_{i+s} - (k_{i+1} - 2 + e)$.  By the previous lemma $\boldsymbol{\mu}$ lies in the lattice 
$$ \Lambda = \Z (p, 0, \dots, 0, -1) \oplus \Z (-1, p, 0, \dots, 0) \oplus \cdots \oplus \Z (0, \dots, 0, -1, p).$$

By \cite{thesispreprint}, Corollary 3.21, we may assume that $w_\tau = 0$ for all $\tau$.
For $j \in \Z / s\Z$, let $c_j = k_j - 2 + p(k_{j-1} - 2) + \cdots + p^{s-1}(k_{j+1} -2 )$.  Assume first that $\theta$ is non-trivial; then $0 < c_j < p^s - 1$.
Let $H$ be an $\F_{p^{2s}}$-vector space scheme over $D^\prime$ defined just as in \cite{thesispreprint}; it satisfies the condition $(\ast \ast )$ of \cite{Ra}.  Let $a_i$, $a_i^\prime$, and $b_i$, for $i \in \Z / 2s\Z$, be parameters defined as in \cite{E}, \S 5 or \cite{thesispreprint}, 4.1.  The relevant facts about them are that $0 \leq a_i^\prime \leq e(p^s - 1)$, that $b_i \in \{ c_{\pi(i)}, 0 \}$ (just as in \cite{thesispreprint}, Lemma 5.3), and that they satisfy the relation
\begin{equation} \label{fundeq}
a_i^\prime = b_{i+1} - pb_i + (p^s - 1) a_i.
\end{equation}

We apply this relation to determine the $a_i$.  As in section 5.1 of \cite{thesispreprint}, we consider four cases:

{\bf {Case 1.}} $b_i = 0, b_{i+1} = c_{i + 1}$.  Then by (\ref{fundeq}) we have 
$$ a_i^\prime - (p^s - 1) a_i = b_{i+1} - pb_i = c_{i+1}.$$
By virtue of the bound on $a_i^\prime$, this equation admits $e$ solutions:
\begin{eqnarray*}
a_i^\prime = c_{i+1} & & a_i = 0 \\
a_i^\prime = c_{i+1} + p^s - 1 & & a_i = 1 \\
\dots & & \dots \\
a_i^\prime = c_{i+1} + (e-1)(p^s - 1) & & a_i = e - 1
\end{eqnarray*}

{\bf {Case 2.}} $b_i = c_i, b_{i+1} = 0$.  Then (\ref{fundeq}) says that
$$ a_i^\prime - (p^s - 1)a_i = -pc_i = \beta - (p^s - 1)(k_{i+1} - 1), $$
where $\beta = (p+1 - k_{i+1}) + p(p+1 - k_i) + \cdots + p^{s-1}(p+1 - k_{i+2})$.  Since $0 < \beta < p^s -1$, we again have $e$ solutions:
\begin{eqnarray*}
a_i^\prime = \beta & & a_i = k_{i+1} - 1 \\
a_i^\prime = \beta + p^s -1 & & a_i = k_{i+1} \\
\dots & & \dots \\
a_i^\prime = \beta + (e-1)(p^s - 1) & & a_i = k_{i+1} - 1 + (e-1)
\end{eqnarray*}

{\bf {Case 3.}} $b_i = 0, b_{i+1} = 0$.  Then $a_i^\prime - (p^s -1)a_i  = 0$, which has $e+1$ solutions:
\begin{eqnarray*}
a_i^\prime = 0 & & a_i = 0 \\
a_i^\prime = p^s - 1 & & a_i = 1 \\
\dots & & \dots \\
a_i^\prime = e(p^s - 1) & & a_i = e
\end{eqnarray*}

{\bf {Case 4.}} $b_i = c_i, b_{i+1} = c_{i+1}$.  Then $a_i^\prime - (p^s - 1)a_i = c_{i+1} - pc_i = -(p^s -1)(k_{i+1} - 2)$, and there are $e+1$ solutions:
\begin{eqnarray*}
a_i^\prime = 0 & & a_i = k_{i+1} - 2 \\
a_i^\prime = p^s - 1 & & a_i = k_{i+1} - 1 \\
\dots & & \dots \\
a_i^\prime = e(p^s - 1) & & a_i = k_{i+1} - 2 + e
\end{eqnarray*}

\begin{lemma}
We may assume without loss of generality that $\{ b_i, b_{i+s} \} = \{ 0, c_i\}$ for each $i \in \Z /2s\Z$.
\end{lemma}
\begin{proof}
We sketch the proof, using the notions and notations of \cite{thesispreprint} without comment.
Recall that $H \subset \Pic^0(\mathbf{M}_{U_1(\p),U}^{bal})[p^\infty]$, where $U \subset G(\A^{\infty, \p})$ is an appropriate open compact subgroup and $\mathbf{M}_{U_1(\p),U}^{bal} \to \Spec D^\prime$ is the semistable model of a Shimura curve as described there and in \cite{Gee}, Thm. 2.18.
 As in \cite{Gee}, $\mathbf{M}_{U_1(\p),U}^{bal}$ represents the functor that associates to an $\mathbf{L}_{1,U}^\ast$-scheme $S$ the collection of canonical balanced $U_1(\p)$-structures on $S$.  The scheme $\mathbf{M}_{U_1(\p),U}^{bal}$ carries an ``Atkin-Lehner'' automorphism $w$ that sends a canonical balanced $U_1(\p)$-structure $(P, P^\prime, \mathcal{K}, \mathcal{K}^\prime)$ to a structure $(Q, Q^\prime, \mathcal{L}, \mathcal{L}^\prime)$, where $\mathcal{L}$ is a lifting of $\mathcal{K}^\prime$ to $\mathbf{E}_{1,U} |_S$ and $Q^\prime$ is the image of $P$ in $\mathcal{L}^\prime$.  The map $w$ interchanges the two components $I$ and $E$ of the special fiber of $\mathbf{M}_{U_1(\p),U}^{bal}$.

By the arguments of \cite{C} \S 10 we see that $\Frob_\p$ preserves $H \oplus w(H)$.  Hence $w(H)$ is an $\F_{p^{2s}}$-vector space scheme over $D^\prime$ lifting the vector space scheme $H_{\phi^q}$ over $K$ on which $\Gal(\overline{K}/K)$ acts via the character $\phi^q$.  Let $w(H)$ be defined by the parameters $a_i^w$, $(a_i^\prime)^w$, $b_i^w$.  Then $a_i^w = a_{i+s}$ and as in \cite{thesispreprint}, Lemma 5.3, we see that $b_i^w = 0$ (resp. $b_i^w = c_i$) if $b_i = c_i$ (resp. $b_i = 0$).

Now, in all the subscripts of the parameters defining $w(H)$, replace $i$ by $i+s$.  We get an $\F_{p^{2s}}$-vector space scheme $\tilde{H}$, defined by parameters $\tilde{a}_i$, $\tilde{a}_i^\prime$, $\tilde{b}_i$, where $\tilde{a}_i = a_i$ and 
$$ \tilde{b}_i = \begin{cases} c_i &: b_{i+s} = 0 \\
0 &: b_{i+s} = c_i \end{cases}$$

Let $N^+ \subset \Z /2s\Z$ (resp. $N^-$) be the set of $i$ such that $b_i = b_{i+s} = c_i$ (resp. $b_i = b_{i+s} = 0$), and let $N = N^+ \cap N^-$.  Suppose first that $N \neq \Z /2s\Z$.  Then there exists an $i$ such that $i \in N$ but $i + 1 \not\in N$.  Suppose that $i \in N^-$ (the case $i \in N^+$ is very similar), and let $n \geq 0$ be the largest integer such that $i - n^\prime \in N^-$ for all $0 \leq n^\prime \leq n$.  
%
It is easy to see that if $\alpha \in \{ i, i+s \}$ is such that $b_{\alpha - n - 1} = c_{\alpha - n - 1}$, then we can switch $b_{i - n^\prime}$ to $c_{i - n^\prime}$ and still obtain the same set of $a_i$'s as possible solutions.  Note that the existence of $\tilde{H}$ guarantees that the $a_{i - n^\prime}$ are in the range where this is possible.  Iterating this procedure proves the lemma.

Finally suppose that $N = \Z / 2s\Z$.  Since $\phi$ is a character of level $2s$, there is some $i$ such that $a_i \neq a_{i+s}$.  We leave it as an exercise to the reader to show that, after possible replacing $i$ with $i+s$, for all $0 \leq n^\prime \leq s - 1$, if $b_{i - n^\prime} = 0$ (resp. $b_{i - n^\prime} = c_{i - n^\prime}$) we may change it to $c_{i - n^\prime}$ (resp. to $0$), and still obtain the same set of $a_i$'s as possible solutions.
\end{proof}

From the definition of $\boldsymbol{\mu}$ we see that $-e \leq \mu_i \leq e$ for all $i$ and that for some $i$ we have $-(e-1) \leq \mu_i \leq e - 1$.  Since $\boldsymbol{\mu} \in \Lambda$, this implies $\boldsymbol{\mu} = 0$.  Thus $a_i + a_{i+s} = k_{i+1} - 2 + e$ for all $i$.

\begin{proposition} \label{mainprop}
Let $\rho: \Gal(\overline{F} / F) \to \GL_2(\Fpbar)$ be such that $\rho |_{G_\p}$ is irreducible and $\rho$ is modular of weight $\sigma$ such that $\sigma_\p$ is a constituent of $Ind_{B}^{\GL_2(k_\p)} \theta$, where $\theta: B \to \Fpbar$ is non-trivial and has the form of (\ref{deftheta}) above.  Then there exists a subset $S \subset I$ and a labeling $\{ \tilde{\tau}, \tilde{\tau}^\prime \}$ of the two liftings of $\tau: k_\p \hookrightarrow \Fpbar$ to $\F_{p^{2s}}$ for each $\tau$, such that 
$$ \rho |_{I_{t,\p}} \sim \left( \begin{array}{cc} \phi & 0 \\ 0 & \phi^q \end{array} \right) , $$
where for each $\tau: k_\p \hookrightarrow \Fpbar$ there is an integer $0 \leq \delta_\tau \leq e - 1$ such that
$$ \phi = \prod_{\tau \in I} (\psi_{\tilde{\tau}} \psi_{\tilde{\tau}^\prime})^{w_\tau} \prod_{\tau \in S} \psi_{\tilde{\tau}}^{k_\tau - 2 + \delta_\tau + \nu_S(\tau)} \psi_{\tilde{\tau}^\prime}^{e-1-\delta_\tau} \prod_{\tau \not\in S} \psi_{\tilde{\tau}}^{p + e - 1 - \delta_\tau} \psi_{\tilde{\tau}^\prime}^{k_\tau - 2 + \delta_\tau + \nu_S(\tau)}.$$
\end{proposition}
\begin{proof}
This is analogous to Proposition 5.6 and Corollary 5.8 of \cite{thesispreprint}.  As in that paper, we reduce to the case of $w_\tau = 0$ for all $\tau \in I$.  Let $\Phi(\theta)$ be the set of all $\phi$ of the form in the statement.  Any $\phi \in \Phi(\theta)$ is specified by the data $(S, \epsi_j, \delta_j)$, where $S \subset I$ and for any $j \in \Z / j\Z$ we have a bijection of two-element sets $\epsi_j: \pi^{-1}(j) = \{ j, j+s \} \to \{ \psi_{\tilde{\tau}_j}, \psi_{\tilde{\tau}_j^\prime} \}$ and an integer $0 \leq \delta_j \leq e-1$.  The character corresponding to $(S, \epsi_j, \delta_j)$ is $\phi = \prod_{i \in \Z /2s\Z} \psi_i^{m_i}$, where
$$ m_i = \begin{cases}
k_i -2 + \nu_S(\tau_i) + \delta_i &: \tau_i \in S, \epsi_i(i) = \psi_{\tilde{\tau}_i} \\
e-1 - \delta_i &: \tau_i \in S, \epsi_i(i) = \psi_{\tilde{\tau}_i^\prime} \\
p + e - 1 - \delta_i &: \tau_i \not\in S, \epsi_i(i) = \psi_{\tilde{\tau}_i} \\
k_i - 2 + \nu_S(\tau_i) + \delta_i &: \tau_i \not\in S, \epsi_i(i) = \psi_{\tilde{\tau}_i^\prime}
\end{cases}.
$$
Here we make the usual abuse of notation: $\tau_i = \tau_{\pi(i)}$, $\delta_i = \delta_{\pi(i)} = \delta_{\tau_i}$, etc.  Clearly every $\phi \in \Phi(\theta)$ is described in this way, although possibly not uniquely. 

Let $\Omega_e(\theta)$ be the set of all $\phi$ satisfying all the conditions emerging from the computations earlier in this section.  Any $\phi \in \Omega_e(\theta)$ is specified by the data $(S^\prime, r_j, \delta_j^\prime)$, where $S^\prime \subset I$ and for every $j \in \Z /s\Z$ we have a bijection $r_j: \{j, j+s \} \to \{ 0, c_j \}$ and an integer $0 \leq \delta_j^\prime \leq e - 1$.  The corresponding character is $\phi = \prod_{i \in \Z /2s\Z} \psi_i^{a_{i-1}}$, where
$$ a_{i-1} = \begin{cases}
e-1 - \delta_i^\prime &: r_i(i) = 0, r_{i+1}(i+1) = c_i \\
k_i - 1 + \delta_i^\prime &: r_i(i) = c_{i-1}, r_{i+1}(i+1) = 0 \\
e-1 - \delta_i^\prime &: r_i(i) = r_{i+1}(i+1) = 0, \tau_{i+1} \in S^\prime \\
e - \delta_i^\prime &: r_{i}(i) = r_{i+1}(i+1) = 0, \tau_{i+1} \not\in S^\prime \\
k_i - 1 + \delta_i^\prime &: r_i(i) = c_{i-1}, r_{i+1}(i+1) = c_i, \tau_{i+1} \in S^\prime \\
k_i - 2 + \delta_i^\prime &: r_i(i) = c_{i-1}, r_{i+1}(i+1) = c_i, \tau_{i+1} \not\in S^\prime
\end{cases}.
$$
Again it is easy to see that every $\phi \in \Omega_e(\theta)$ is described (non-uniquely) in this way.  Here $r_i(i) = 0$ and $r_i(i) = c_{i-1}$ correspond to $b_{i-1} = 0$ and $b_{i-1} = c_{i-1}$, respectively, and $S^\prime$ accounts for the extra possibilities in Cases 3 and 4.  As in \cite{thesispreprint}, Prop.~ 5.6 one constructs a bijection between these two collections of data and deduces that $\Phi(\theta) = \Omega_e(\theta)$.
\end{proof}

\begin{theorem} \label{mainresult}
Suppose that $e < p - 1$ and let $\rho: \Gal(\overline{F} / F) \to \GL_2(\Fpbar)$ be such that $\rho |_{G_\p}$ is irreducible and $\rho$ is modular of weight $\sigma$, where $\sigma_\p$, written as in (\ref{sigmap}), satisfies $k_\tau - 2 + e \leq p - 1$ for all $\tau$.  Then there exists a labeling $\{ \tilde{\tau}, \tilde{\tau}^\prime \}$ of the two liftings of $\tau: k_\p \hookrightarrow \Fpbar$ to $\F_{p^{2s}}$ for each $\tau$, such that 
$$ \rho |_{I_{t,\p}} \sim \left( \begin{array}{cc} \phi & 0 \\ 0 & \phi^q \end{array} \right) , $$
where for each $\tau: k_\p \hookrightarrow \Fpbar$ there is an integer $0 \leq \delta_\tau \leq e - 1$ such that
$$ \phi = \prod_{\tau \in I} (\psi_{\tilde{\tau}} \psi_{\tilde{\tau}^\prime})^{w_\tau} \prod_{\tau \in I} \psi_{\tilde{\tau}}^{k_\tau - 1 + \delta_\tau} \psi_{\tilde{\tau}^\prime}^{e-1-\delta_\tau} .$$
\end{theorem}
\begin{proof}
As in \cite{thesispreprint} we may assume that $w_\tau = 0$ for all $\tau$.  Assume first that $k_\tau \neq 2$ for some $\tau$.  Denote by $\Theta(\sigma_\p)$ the set of all characters $\theta: B \to \Fpbar^\ast$ such that $\sigma_\p$ is a constituent of $Ind_B^{\GL_2(k_\p)} \theta$; all these characters $\theta$ are non-trivial.  The elements of $\Theta(\sigma_\p)$ are the following, where $T$ runs over all $T \subset I$:
$$\theta_T : \left( \begin{array}{cc} a & b \\ 0 & d \end{array} \right) \mapsto \prod_{\tau \in T} \tau(ad)^{p-1} \tau(d)^{k_\tau - 1 - \nu_T(\tau)} \prod_{\tau \not\in T} \tau(ad)^{k_\tau - 2} \tau(d)^{p + 1 - k_\tau - \nu_T(\tau)}.$$

If $\rho$ is modular of a weight whose $\p$-component is $\sigma_\p$, then $\phi \in \bigcap_{\theta \in \Theta(\sigma_\p)} \Phi(\theta)$, and we will compute this intersection.  If $s = 1$, then the desired result is immediate from Proposition \ref{mainprop} by considering $\Phi(\theta_I)$.  Otherwise, suppose that $\phi \in \bigcap_{\theta \in \Theta(\sigma_\p)} \Phi(\theta)$, but $\phi$ is not of the form specified in the statement of the theorem.  Since $\phi \in \Phi(\theta_I)$, it is easy to see that $\phi = \prod_{j \in \Z /2s\Z} \psi_j^{m_j}$ where $\{ m_i, m_{i+s} \} = \{ \epsi_i, k_i - 2 + e - \epsi_i \}$, where $0 \leq \epsi_i \leq e$ and for some $i$ we have $\epsi_i = e$.  Moreover, we may assume that $k_i > e + 1$, since otherwise $\{ k_i - 2, e \} = \{ k_i - 2 + e - \epsi_i, \epsi_i \}$ for some $0 \leq \epsi_i \leq e - 1$.

If $s \geq 2$, then the elements of $\Phi(\theta_{T = \{ \tau_i \} })$ are the following, as $S$ runs over the subsets of $I$ and each $\delta_\tau$ runs over $\{ 0, 1, \dots, e - 1 \}$:
\begin{multline*}
\prod_{\tau \in S \atop \tau \neq \tau_i} \psi_{\tilde{\tau}}^{k_\tau - 2 + e - \delta_\tau - \nu_S(\tau)} \psi_{\tilde{\tau}^\prime}^{p + \delta_\tau - \nu_T(\tau)}
\prod_{\tau \not\in S \atop \tau \neq \tau_i} \psi_{\tilde{\tau}}^{k_\tau - 2 + e - \delta_\tau - \nu_S(\tau)} \psi_{\tilde{\tau}^\prime}^{\delta_\tau - \nu_T(\tau)} \\
\times \begin{cases}
\psi_{\tilde{\tau}_i}^{k_i + p - 1 + \delta_i} \psi_{\tilde{\tau}_i^\prime}^{p + e - 1 - \delta_i - \nu_S(\tau_i)} &: \tau_i \in S \\
\psi_{\tilde{\tau}_i}^{k_i - 1 + \delta_i} \psi_{\tilde{\tau}_i^\prime}^{p + e - 1 - \delta_i - \nu_S(\tau_i)} &: \tau_i \not\in S
\end{cases}
\end{multline*}

Dividing this by the expression for $\phi$ found above, we see that for some $S \subset I$ we have
\begin{multline*}
1 = 
\prod_{\tau \in S \atop \tau \neq \tau_i} \psi_{\tilde{\tau}}^{\epsi_\tau - \delta_\tau - \nu_S(\tau) \atop k_\tau - 2 + e - \delta_\tau - \epsi_\tau - \nu_S(\tau)} \psi_{\tilde{\tau}^\prime}^{ p + \delta_\tau - \epsi_\tau - \nu_T(\tau) \atop  p + \delta_\tau - k_\tau + 2 - e + \epsi_\tau}
\prod_{\tau \not\in S \atop \tau \neq \tau_i} \psi_{\tilde{\tau}}^{\epsi_\tau - \delta_\tau - \nu_S(\tau) \atop k_\tau - 2 + e - \delta_\tau - \epsi_\tau - \nu_S(\tau)} \psi_{\tilde{\tau}^\prime}^{\delta_\tau - \epsi_\tau - \nu_T(\tau) \atop  \delta_\tau - k_\tau + 2 - e + \epsi_\tau} \\
\times \begin{cases}
\psi_{\tilde{\tau}_i}^{p + 1 + \delta_i \atop p + k_i - 1 - e + \delta_i} \psi_{\tilde{\tau}_i^\prime}^{p - 1 - \delta_i - \nu_S(\tau_i) \atop p + 1 + e - k_i - \delta_i - \nu_S(\tau_i)} &: \tau_i \in S \\
\psi_{\tilde{\tau}_i}^{1 + \delta_i \atop k_i - 1 - e + \delta_i} \psi_{\tilde{\tau}_i^\prime}^{p - 1 - \delta_i - \nu_S(\tau_i) \atop p + 1 + e - k_i - \delta_i - \nu_S(\tau_i)} &: \tau_i \not\in S
\end{cases}
\end{multline*}
Here for each pair $\psi_{\tilde{\tau}}, \psi_{\tilde{\tau}^\prime}$ we choose either the top or the bottom exponent in both cases.  If we rewrite this expression as $\prod_{i \in \Z /2s\Z} \psi_i^{r_i}$, then we must have $(r_0, \dots, r_{2s-1}) \in \Lambda$.  Under our hypotheses, all these exponents lie in the range $[ - (p-1), 2p - 2]$.  However, they cannot all be $-(p-1)$, nor can they all be $2p - 2$, and hence the only possible values of the $r_i$ are $-1, 0, p-1$, and $p$.  Consider now the exponent $r_{\tilde{\tau}_i}$ of $\psi_{\tilde{\tau}_i}$.  Since $1 \leq 1 + \delta_i \leq p - 2$ and $1 \leq k_i - 1 - e + \delta_i \leq  p - 2$ (recall $k_i > e + 1$), we see that $r_{\tilde{\tau}_i}$ cannot take any of the allowed values, whence we cannot have $\tau_i \not\in S$.  But similarly $\tau_i \in S$ is impossible.  We obtain a contradiction, which proves that $\phi \not\in \bigcap_{\theta \in \Theta(\sigma_\p)} \Phi(\theta)$.

Finally, suppose $k_\tau = 2$ for all $\tau$.  In this case (recall $w_\tau = 0$ for all $\tau$) the only $\theta$ such that $\sigma_\p$ is a constituent of $Ind_B^{\GL_2(k_\p)} \theta$ is the trivial character (\cite{Fredpreprint}, Prop. 1.1).  Just as in \cite{thesispreprint}, 5.4., we construct an $\F_{p^{2s}}$-vector space scheme $V$ such that $\Gal(\overline{K}/K)$ acts on $V_K$ by the character $\phi$.  Let $a_i, a_i^\prime, b_i$ be the parameters associated to $V$.  As in \cite{thesispreprint} we see that $b_i = 0$ for all $i$; hence, by (\ref{fundeq}), each $a_i$ can take any value between $0$ and $e$.  Our claim now follows from Lemma \ref{centchar}.
\end{proof}
\begin{remark}
As in \cite{thesispreprint} \S 5, it is possible to relax the hypothesis that $k_\tau - 2 + e \leq p - 1$ for all $\tau$, at the price of obtaining a somewhat weaker result.  In this case, the set $\bigcap_{\theta \in \Theta(\sigma_\p)} \Phi(\theta)$ will be larger than the conjectured set of $\phi$'s for representations modular of a weight with $\p$-component $\sigma_\p$.  However, we can still assert that $\phi \in \bigcap_{\theta \in \Theta(\sigma_\p)} \Phi(\theta)$.
\end{remark}

\section{Examples}
Let $F = \Q(\sqrt{5})$.  Let $p = 5$; then $(p) = \p^2$ in $F$, where $\p = ((5 + \sqrt{5})/2)$, and $k_\p = \F_5$.  Thus we have $e = 2$ and $s = 1$.  The weights in this situation are $ \det^{w} \Sym^{k-2} \F_5 \tensor \overline{\F}_5 = F(w + k - 2, w)$, where $2 \leq k \leq 6$ and $0 \leq w \leq 3$.  All our examples rely on Lassina Demb\'{e}l\'{e}'s computations of Hilbert modular forms (see \cite{Dembele}), which so far exist only for $\Q(\sqrt{5})$.  For each Hilbert modular form, Demb\'{e}l\'{e} computes the list of weights for which the associated mod $5$ Galois representation $\overline{\rho}$ is modular.  He also provides evidence for (but does not actually compute) the projective image of $\overline{\rho}^{ss}$; clearly $\overline{\rho}$ is reducible if and only if this projective image is cyclic.

We have used Magma to find (elliptic) modular newforms $f$ with integer coefficients.  Then $\overline{\rho}_{f} |_{I_p}$ is described by classical theorems of Deligne and Fontaine.  We search for the base change of $f$ to $F$ in Demb\'{e}l\'{e}'s tables and obtain the weights for which $\overline{\rho}_f |_{\Gal(\overline{F}/F)}$ is modular.  In all examples that we have computed the results are, fortunately, consistent with Conjecture \ref{mainconj}.

\subsection{Non-ordinary forms}
If $f$ is non-ordinary at $5$ and has weight $2 \leq k \leq 6$, then $\overline{\rho}_f$ is tame at $5$.  By a result of Fontaine (see \cite{E}, Thm. 2.6),
$$ \overline{\rho} |_{I_5} \sim \left( \begin{array}{cc} \psi^{k - 1} & 0 \\ 0 & \psi^{5(k-1)} \end{array} \right), $$
where $\psi$ is a fundamental character of level $2$.  From the description of the isomorphism between $I_{t,\p}$ and $\varprojlim \F_{p^n}^\ast$ (see, for instance, \cite{thesispreprint}, 4.1) we see that
$$ \overline{\rho} |_{I_\p} \sim \left( \begin{array}{cc} \psi^{2(k - 1)} & 0 \\ 0 & \psi^{10(k-1)} \end{array} \right). $$
The weights predicted by our conjecture are the following:
\begin{eqnarray*}
F(0,0), F(3,1), F(4,0), F(5,3), & & k = 2,6 \\
F(1,1), F(4,2), F(5,1), F(6,4), & & k = 5 \\
F(2,0), F(3,3), F(4,2), F(7,3), & & k = 3 \\
F(0,0), F(3,1), F(4,0), & & k = 4 
\end{eqnarray*}

Here are some of the computational results.  Observe that the form with $k = 4$ gives a tame example of level $1$, where the associated local Galois representation at $\p$ is scalar.  In this case the global mod $5$ Galois representation is reducible; hence this is not a counterexample to the conjecture, even though we obtain only two weights.  The other representations in the list are irreducible.

\vskip0.1cm
\begin{tabular}{|r | r | l | | l|} \hline
weight & level & $q$-expansion of $f$ & modular weights of $\overline{\rho}_f |_{\Gal(\overline{F}/F)}$ \\
\hline\hline
$2$ & $14$ & $q - q^2 - 2q^3 + q^4 + 2q^6 + q^7 + O(q^8)$ & $F(0,0), F(3,1), F(5,3), F(4,0)$ \\ \hline
$3$ & $7$ & $q - 3q^2 + 5q^4 - 7q^7 + O(q^8)$ & $F(3,3), F(2,0), F(4,2), F(7,3)$ \\ \hline
$3$ & $8$ & $q - 2q^2 - 2q^3 + 4q^4 + 4q^6 + O(q^8)$& $F(3,3), F(2,0), F(4,2), F(7,3)$ \\ \hline
$4$ & $9$ & $q - 8q^4 + 20q^7 - 70q^{13} + O(q^{16})$ & $F(0,0), F(4,0)$ \\ \hline
$6$ & $14$ & $q + 4q^2 + 8q^3 + 16q^4 + 10q^5 + 32q^6 + O(q^7)$ & $F(0,0), F(3,1), F(5,3), F(4,0)$ \\ \hline
\end{tabular}

\subsection{Ordinary forms}
Elliptic modular newforms which are ordinary at $5$ are much more plentiful than non-ordinary ones.  In this case, $\overline{\rho}_f |_{\Gal(\overline{F}/F)}$ is not in general tame at $\p$, and it is natural to expect that even when $\overline{\rho}_f |_{\Gal(\overline{F}/F)}$ is irreducible, the modular weights will be only a subset of those which are modular for the semisimplification.  If $f$ has weight $2 \leq k \leq 6$, then by a theorem of Deligne (see \cite{E}, Thm. 2.5),
$$ \overline{\rho}_f |_{I_\p} \sim \left( \begin{array}{cc} \psi^{2(k-1)} & \ast \\ 0 & 1 \end{array} \right),$$
where $\psi$ is a fundamental character of level $1$.  If $\overline{\rho}_f$ is tame, then Conjecture \ref{mainconj} predicts the following sets of weights:
\begin{eqnarray*}
F(0,0), F(2,2), F(3,1), F(5,3), F(4,0), F(6,2),& & k = 2,4,6 \\
F(3,3), F(2,0), F(7,3),& & k = 3,5
\end{eqnarray*}

For most of the forms we found, Demb\'{e}l\'{e}'s computations suggest that the global Galois representation is reducible.  We found only three irreducible examples, which are compatible with the conjecture:
\vskip0.1cm
\begin{tabular}{|r | r | l | | l|} \hline
wt. & level & $q$-expansion of $f$ & modular weights \\
\hline\hline
$4$ & $8$ & $q - 4q^3 - 2q^5 + 24q^7 + O(q^9)$ & $F(0,0), F(4,0)$ \\ \hline
$6$ & $8$ & $q + 20q^3 - 74q^5 + O(q^7)$ & $F(4,0)$ \\ \hline
$6$ & $9$ & $q + 6q^2 + 4q^4 - 6q^5 + O(q^7)$ & $F(4,0)$ \\ \hline
\end{tabular}
\vskip0.1cm

In the examples where the global Galois representation appears to be reducible, we can still apply our conjecture to $\overline{\rho}_f |_{I_\p}$ to obtain a set $W_\p^?(\overline{\rho}_f)$.  The computed modular weights always lie inside this set.
\vskip0.1cm
\begin{tabular}{|r | r | l | | l|} \hline
wt. & level & $q$-expansion of $f$ & modular weights \\
\hline\hline
$2$ & $8$ & $q + 2q^2 + 2q^3 + 4q^4 + 4q^5 + 4q^6 + O(q^7)$ & $F(3,1), F(5,3)$ \\ \hline
$2$ & $9$ & $q - 3q^2 + 7q^4 - 6q^5 + O(q^7)$ & $F(3,1), F(5,3)$ \\ \hline
$3$ & $3$ & $q + 3q^2 + 9q^3 + 13q^4 + 24q^5 + 27q^6 + O(q^7)$ & $F(3,3), F(7,3)$ \\ \hline

$3$ & $4$ & $q + 4q^2 + 8q^3 + 16q^4 + 26q^5 + 32q^6 + O(q^7)$ & $F(3,3), F(2,0), F(7,3)$ \\ \hline
$3$ & $7$ & $q + 5q^2 + 8q^3 + 21q^4 + 24q^5 + 40q^6 + O(q^7)$ & $F(3,3), F(7,3)$ \\ \hline
$3$ & $8$ & $q + 4q^2 + 10q^3 + 16q^4 + 24q^5 + 40q^6 + O(q^7)$ & $F(3,3), F(7,3)$ \\ \hline
$4$ & $6$ & $q - 2q^2 - 3q^3 + 4q^4 + 6q^5 + 6q^6 + O(q^7)$ & $F(2,2), F(6,2)$ \\ \hline
$4$ & $8$ & $q - 4q^3 - 2q^5 + 24q^7 + O(q^9)$ & $F(0,0), F(4,0)$ \\ \hline
$4$ & $8$ & $q + 8q^2 + 26q^3 + 64q^4 + 124q^5 + 208q^6 + O(q^7)$ & $F(3,1), F(5,3)$ \\ \hline
$4$ & $9$ & $q - 9q^2 + 73q^4 - 126q^5 + O(q^7)$ & $F(3,1), F(5,3)$ \\ \hline
$5$ & $4$ & $q - 4q^2 + 16q^4 - 14q^5 + O(q^8)$ & $F(3,3), F(2,0), F(7,3)$ \\ \hline
$5$ & $7$ & $q + 17q^2 + 80q^3 + 273q^4 + 624q^5 + 1360q^6 + O(q^7)$ & $F(3,3), F(7,3)$ \\ \hline
$5$ & $8$ & $q + 16q^2 + 82q^3 + 256q^4 + 624q^5 + 1312q^6 + O(q^7)$ & $F(3,3), F(7,3)$ \\ \hline
$6$ & $8$ & $q + 20q^3 - 74q^5 + O(q^7)$ & $F(4,0)$ \\ \hline
$6$ & $8$ & $q + 32q^2 + 242q^3 + 1024q^4 + 3124q^5 + 7744q^6 + O(q^7)$ & $F(3,1), F(5,3)$ \\ \hline
$6$ & $9$ & $q + 6q^2 + 4q^4 - 6q^5 + O(q^7)$ & $F(4,0)$ \\ \hline
$6$ & $9$ & $q - 33q^2 + 1057q^4 - 3126q^5 + O(q^7)$ & $F(3,1), F(5,3)$ \\ \hline
\end{tabular}

\bibliographystyle{math}
\bibliography{wt4}

\end{document}